\documentclass[number,citesort,dvips]{arxbj}
\usepackage{upgreek}

\aid{0}
\volume{17}
\issue{2}
\pubyear{2011}
\firstpage{545}
\lastpage{561}
\doi{10.3150/10-BEJ283}

\makeatletter

\newremark{Example}{Example}
\newtheorem{theorem}{Theorem}
\newtheorem{lemma}{Lemma}
\newremark{remark}{Remark}

\renewcommand{\triangle}{\Delta}
\newcommand{\cdott}{}
\makeatother

\begin{document}
\begin{frontmatter}

\title{Central limit theorems for local empirical processes near boundaries of sets}
\runtitle{Empiricals near set boundaries}

\begin{aug}
\author[a]{\fnms{John H.J.} \snm{Einmahl}\thanksref{a}\ead[label=e1]{j.h.j.einmahl@uvt.nl}}
\and
\author[b]{\fnms{Est\'ate V.} \snm{Khmaladze}\thanksref{b}\ead[label=e2]{estate.khmaladze@vuw.ac.nz}}

\runauthor{J.H.J. Einmahl and E.V. Khmaladze}

\address[a]{Department of Econometrics and OR and CentER,
Tilburg University,
P.O. Box 90153,
5000 LE Tilburg,
The Netherlands.
\printead{e1}}

\address[b]{School of Mathematics,
Statistics and Operations Research,
Victoria University of Wellington,
P.O. Box 600,
Wellington,
New Zealand.
\printead{e2}}
\end{aug}

\received{\smonth{2} \syear{2009}}
\revised{\smonth{7} \syear{2009}}

%
\begin{abstract}
We define the local empirical process, based on $n$ i.i.d. random
vectors in dimension $d$, in the neighborhood of the boundary of a
fixed set. Under natural conditions on the shrinking neighborhood, we
show that, for these local empirical processes, indexed by classes of
sets that vary with $n$ and satisfy certain conditions, an
appropriately defined uniform central limit theorem holds. The concept
of differentiation of sets in measure is very convenient for
developing the results. Some examples and statistical applications are
also presented.
\end{abstract}

%

\end{frontmatter}

\section{Introduction}\label{sec1}

Let $X_1,\dots, X_n$ be independent and identically distributed
(i.i.d.) random vectors in $\mathbb{R}^d$ ($d \in
\mathbb{N}$), distributed according to an absolutely continuous
probability measure $P$. Denote the corresponding density by $p$. For a
Borel measurable subset $D$ of $\mathbb{R}^d$, write
\[
\Psi_n(D)=\sum_{i=1}^n 1_D (X_i).
\]
The process $\Psi_n(D), D \in\mathcal{D}$ ($\mathcal{D}$ being the
class of Borel
sets), is, by definition, a binomial process on $\mathbb{R}^d$; $\Psi
_n/n$ is the empirical measure corresponding to
$X_1,\dots, X_n$. Clearly, $\mathbb{E}\Psi_n(D)=nP(D)$.

Let $K$ be a convex body in $\mathbb{R}^d$. The set $K$ will be \textit
{fixed} throughout. Denote its boundary by $\partial K$. It is the aim
of this paper to study the behavior of $\Psi_n$ in the neighborhood of
$\partial K$. Write $\Vert z-\partial K \Vert=\min_{x\in
\partial K} \Vert z-x\Vert$ and let
\[
\mathcal{V}_{\varepsilon}(\partial K)=\{z\in\mathbb{R}^d\dvtx
{\Vert
z-\partial K \Vert} \leq\varepsilon\},\qquad \varepsilon>0 ,
\]
denote this neighborhood. Set $a=P(\mathcal{V}_{\varepsilon}(\partial K))$.
For a Borel set
$A\subset\mathcal{V}_{\varepsilon}(\partial K)$, define
\[
z_n(A)=\frac{1}{\sqrt{na}}[\Psi_n(A)-nP(A)].
\]

If $\varepsilon\to0$, then all sets $A$ will ``shrink toward''
$\partial K$.
If, however, $n\to\infty$ at the same time, the random variables
$z_n(A)$ do not have to converge to 0 and, if $n\varepsilon\to\infty
$, they
should typically converge to Gaussian random variables. However, where
would these Gaussian random variables ``live''? Would they form some
set-parametric process? These are the questions we seek to address in
this paper. We will do this using the concept of differentiability of
set-valued functions, as was recently developed in \cite{K07}.

One cannot prove a sufficiently interesting Gaussian limit theorem for
$z_n$ which is defined on all Borel subsets of $\mathcal
{V}_{\varepsilon
}(\partial K)$. Instead, one needs to consider smaller classes of sets.
Let $\varepsilon=\varepsilon_n\to0$ as $n\to\infty$ and let
$\mathcal{A}_{\varepsilon
_n}$ be a class of measurable subsets of $\mathcal{V}_{\varepsilon
_n}(\partial K)$. The \textit{canonical example} of $\mathcal
{A}_{\varepsilon
_n}$ is
constructed as follows. Let $\mathcal{K}$ be a fixed class of Borel
sets in $\mathbb{R}^d$ and define $\mathcal{A}=\{K'\triangle K\dvtx  K'
\in\mathcal{K}\}$, where $\triangle$ denotes ``symmetric
difference''. Now, take $\mathcal{A}_{\varepsilon_n}=\{A\in\mathcal
{A}\dvtx  A
\subset\mathcal{V}_{\varepsilon_n}(\partial K)\}$. Particular cases
can be
found in Examples \ref{Example1} and \ref{Example2} below.

Our main result is the central limit theorem for the \textit{local
empirical process near} $\partial K$ \textit{and indexed by} $\mathcal
{A}_{\varepsilon_n}$,
%
%
\begin{equation}\label{lep} \{z_n(A), A\in\mathcal{A}_{\varepsilon
_n}\}.
\end{equation}
Denoting the conditional probability distribution on $\mathcal
{V}_{\varepsilon
}(\partial K)$ by $P_\varepsilon(A)=P(A)/a$, we can also write
\[
z_n(A)=\frac{1}{\sqrt{na}}[\Psi_n(A)-naP_\varepsilon(A)] .
\]
This reflects the fact that, on average, the effective sample size is
equal to $na$, not $n$. We therefore assume, in addition to
$\varepsilon_n\to
0$, that
\[
n\varepsilon_n\to\infty\qquad\mbox{as } n\to\infty.
\]
This will imply that $na\to\infty$ and ensure that the sets in
$\mathcal{A}_{\varepsilon_n}$ contain enough observations to obtain Gaussian
limit behavior.

Although very natural here, it is, in general, unusual that an
empirical process is defined on a class of sets that depends on $n$. We
will show that its limiting process should be defined on a class of
subsets not of the ``same'' $\mathbb{R}^d$, but of the cylinder
$\partial K \times[-1,1]$. The subsets in this class are properly
defined derivatives of sequences of sets, with the $n$th set an element
of $\mathcal{A}_{\varepsilon_n}$.

Poisson limit behavior of $\Psi_n$ on $\mathcal{V}_{\varepsilon
_n}(\partial
K)$ has been studied in \cite{KW08}. The main limit
result there had a somewhat unusual property:
it contained a functional limit theorem, but not a one-dimensional
limit theorem for $\Psi_n$.
Indeed, although it showed weak convergence and, moreover, convergence
in total variation, for the process $\Psi_n$ given on \textit{all} Borel
subsets of $\mathcal{V}_{\varepsilon_n}(\partial K)$, for a particular
sequence of subsets $A_{\varepsilon_n}$, it remained unspecified which random
variable from the limiting process the sequence $\Psi_n(A_{\varepsilon_n})$
would converge to. This happened because the notion of derivative sets
had not been developed at the time Khmaladze and Weil \cite{KW08} was
accepted for publication. In this paper, the situation is different --
extracting the one-dimensional limit theorem from Theorem \ref{theo1} gives the
following statement: \textit{if the set-valued function $A_{\varepsilon
}$ is
differentiable in $\varepsilon$ at $\varepsilon=0$ and
$\mathrm{d}A_{\varepsilon}/\mathrm{d}\varepsilon$ is its
derivative} (see Section \ref{sec3} or \cite{K07}), \textit{then}
\[
z_n(A_{\varepsilon_n})\stackrel{d}{\to} W(\mathrm{d}A_{\varepsilon
}/\mathrm{d}\varepsilon),
\]
where $W$ is the set-parametric Brownian motion defined just before
Lemma \ref{lemma1}, Section~\ref{sec4}.

The local empirical process for one-dimensional $X_i$, that is, the
empirical process in the neighborhood of a point $c \in\mathbb{R}\cup
\{-\infty, \infty\}$, is a classical object in probability theory,
one which has proven to be very valuable in statistics; see, for
example, \cite{CH93,M88,DEH89,DM90,DM91,K98,E92},
the book by Cs\"org\H{o} and Horv\'ath \cite{CH93} and Khmaladze \cite{K98}.
The one-dimensional local
empirical process has been extended to the multivariate set-up, but,
typically, only the neighborhood of a point $c \in\mathbb{R}^d$ or
the region outside a large sphere are considered; see, for example,
\cite{DM94,R94,DH98,M04,DZ08}.
Perhaps the closest to the
present paper
are \cite{EM97} and \cite{E97}. For a local empirical
process for function-valued random elements, see \cite{EL06}.

The paper is organized as follows. In the next section, we present
statistical applications. In Section \ref{sec3}, we introduce the necessary
geometry and the appropriate concept of differentiation of sets. In
Section \ref{sec4}, the main results, central limit theorems for $z_n$, and some
examples will be presented. Proofs are collected in Section \ref{sec5}.

\section{Statistical motivation}\label{sec2}

Although the local empirical process near the boundary of a set is an
interesting probabilistic object in its own right, the study of this
type of process was mostly motivated by
problems in spatial statistics. Consider a family of distributions,
indexed by some parameter $\theta$, and denote by $L_n(\theta, \theta
')$ the log-likelihood ratio. If the parameter were a vector, as in
parametric problems (see, e.g., \cite{IH81}),
the local analysis of $L_n(\theta,\theta')$ (or any other process
which the inference is based upon), in $\theta'$ from the neighborhood
of the true value $\theta$, is a crucial step in asymptotic
statistical theory. It forms, for example, the basis of contiguity
theory. The situation is similar when the parameter is a function (see,
e.g., \cite{BicEtal93,W06}).
However, it has thus far not been known how to carry out such a local
analysis when the parameter is a set.

Examples of set-parametric problems are provided by the class of
spatial change point problems or \textit{change set problems} (see, e.g.,
\cite{KMT06}). In these problems, the
observation is usually a (marked) point process in ${\mathbb R}^d$ and
the model assumption is that there is a set, or an image, $K$, such
that outside $K$, the distribution of the point process (e.g., the
distribution of the marks) sharply changes.
One can think of $K$ as, for example, an ore deposit site, a pollution
site or a site with different magnetic properties. The literature on
this problem is very broad; see, for example,
\cite{KT93,F04,MS99,CFR07}.

In most of the particular formulations of the change set problem, the
log-likelihood ratio $L_n(K,K')$ is some form or another of the local
empirical process (\ref{lep}), where $K$ plays the role of the true
value of the change set, while the sets $K'$ are small deviations from it.
To be more precise, let $(X_1,Y_1), \ldots, (X_n,Y_n)$ be independent
random vectors, with $X_i$ being (as before) a $d$-dimensional location
and $Y_i$ being a ``mark'', not necessarily one-dimensional. Write
$K(\varepsilon)$ instead of $K'$ in order to explicitly express the
dependence on $\varepsilon>0$. Let $P_1$ and $P_2$ be the
distributions of $Y_i$ outside $K$ and on $K$, respectively. The
log-likelihood ratio
then has the form
\[
L_n(K,K(\varepsilon))=\sum_{i=1}^n\bigl[{1}_{K(\varepsilon)\setminus
K}(X_i) - {1}_{K \setminus K(\varepsilon)}(X_i) \bigr] \xi(Y_i) ,
\]
where
$\xi(Y_i)=\log\frac{\mathrm{d}P_2}{\mathrm{d}P_1}(Y_i)$. We focus on
%
%
\begin{equation}\label{plmi}\sum_{i=1}^n{1}_{K(\varepsilon)\setminus
K}(X_i) \quad\mbox{and}\quad \sum_{i=1}^n {1}_{K \setminus
K(\varepsilon)}(X_i) ;
\end{equation}
a discussion of the behavior of the $\xi(Y_i)$ is of secondary
importance here. Let $\tilde{ \mathcal{K}}$
be a class of set-valued functions $K(\cdot)$ all converging to the
same $K$, that is, $K(\varepsilon)\triangle K$ shrinks toward the
boundary $\partial K$ when $\varepsilon\to0$. These sets describe
the deviations from the hypothetical change set $K$. Let us consider
the processes in (\ref{plmi}) given on these deviations and
investigate their joint limit in distribution when $n \to\infty$ and
$\varepsilon=\varepsilon_n\to0$: the larger the number of
observations, the smaller, or narrower, the sets we consider.
In the appropriate formulation of ``local alternatives'', where not only
$K(\varepsilon)$ tends to $K$, but also $P_2$ tends to $P_1$, when $n
\to\infty$, the convergence to a Gaussian limiting process -- as
shown in this paper -- is of fundamental importance.
Indeed, although there is a rich literature on statistical estimation
of sets, we know very few results on testing hypotheses about sets and
no results for testing against local alternatives.

To illustrate another class of statistical problems where the parameter
is a set, consider two prominent examples: the excess mass approach
(cf. \cite{MS91,N91}), and the shorth
\cite{AndEtal72,G88}) and its
generalization \cite{EM92}. Let $\mathcal{K}$ be a
fixed subset of $\mathcal{D}$, as in the canonical example. One could
choose, for instance, $\mathcal{K}$ to be the class of all ellipsoids.
Define the excess mass set $K$ for level $\lambda>0$ by
\[
K=\mathop{\arg\max}_{K'\in\mathcal{K}}\{P(K')-\lambda\mu_d(K')\},
\]
where $\mu_d$ denotes $d$-dimensional Lebesgue measure. Similarly, the
generalized shorth or minimum volume set
$K$ for probability $\alpha\in(0,1)$ is given by
\[
K=\mathop{\arg\min}_{K'\in\mathcal{K}}\{\mu_d(K')\dvtx P(K')\geq
\alpha\};
\]
see \cite{R85,D92} when $\mathcal{K}$ is the class
of all ellipsoids.
It turns out that both of these sets $K$ and their M-estimators can be
analyzed somewhat similarly. Therefore, we confine ourselves to the
excess mass set.

The obvious non-parametric estimator for this set is obtained by
replacing $P$ by $\Psi_n/n$:
\begin{eqnarray*}
K_n&=&\mathop{\arg\max}_{K'\in\mathcal
{K}}\{\Psi_n(K')/n-\lambda\mu_d(K')\}\\
&=&\mathop{\arg\max}_{K'\in\mathcal{K}}n^{2/3}\bigl\{\Psi
_n(K')/n-\lambda\mu_d(K')-\bigl(\Psi_n(K)/n-\lambda\mu_d(K)\bigr)\bigr\}.
\end{eqnarray*}
Let $\varepsilon_n$ be such that $a=n^{-1/3}$ (cf. the ``cube root
asymptotics'' of \cite{KP90}). Under certain conditions, it
can be shown that for large $T>0$, with high probability,
$K_n\triangle K \in\mathcal{A}_{T\varepsilon_n}$. Observe that in
that case,
\[
K_n=\mathop{\arg\max}_{K'\dvtx  K'\triangle K\in\mathcal
{A}_{T\varepsilon_n}
}\bigl\{z_n(K')-z_n(K)+n^{2/3}\bigl[P(K')-P(K)-\lambda\bigl(\mu_d(K')-\mu_d(K)\bigr)\bigr]\bigr\}.
\]
Now, a central limit theorem for $z_n$ makes it possible to show that,
asymptotically, $K_n\triangle K$
can be described as a deterministic function depending on $n$
(actually, on $\tau^{-1}_{\varepsilon_n}$; see Section \ref{sec4}) evaluated
at a random variable that does not depend on $n$.
This random variable is the $\arg\max$ of some Brownian motion with drift.
Such a result is very useful for a refined analysis of $K_n$. See
\cite{BE10} for a study
of the behavior of such $K_n$'s along these lines.

\section{Some geometry and differentiability of sets}\label{sec3}

In this section, we first briefly review some relevant notation and
facts from geometry. We then recall the concept of ``differentiation of
sets in measure'', as given in \cite{K07}. In that paper and the
references therein (in particular \cite{S93}), more details about
the required geometry can be found. We also refer to the recent
monograph \cite{SW08}.

Let $K\in\mathcal{D}$ be our convex body, that is, a closed, bounded
convex set that has interior points. Denote by $\Pi(z)$ the metric
projection of $z \in\mathbb{R}^d$ on $\partial K$, that is, $\Pi
(z)$ is a nearest point to $z$ on $\partial K$. The set of $z$-values
for which such a nearest point is not unique is a subset $S_K$ of $K$
called the \textit{skeleton} of $K$. Let $\mu_d$ denote
$d$-dimensional Lebesgue measure. It is then known that $\mu_d(S_K)=0$.
A unit vector $u$ is called an \textit{outer normal} of $K$ at $x \in
\partial K$ if there is some $z \in\mathbb{R}^d \backslash K$ such
that $x=\Pi(z)$ and $u=(z-\Pi(z))/\Vert z-\Pi(z) \Vert$.
Let $B_r(z)$ denote the closed ball with center~$z$ and radius $r$. For
$x \in\partial K$, we define the local interior reach
\[
r(x)= \max\{r\dvtx  x \in B_r(z) \subset K\}.
\]
If $r(x)>0$, then the outer normal $u$ at $x \in\partial K$ is unique.
In this case, the unit vector $-u$ is the unique inner normal. In
general, at each $x \in\partial K$, we denote the set of outer normals
by $N(x)$ and the normal bundle of $K$ is defined as
\[
\mathit{Nor}(K)=\{(x,u)\dvtx  x \in\partial K, u \in N(x)\}.
\]
The cylinder $\Sigma=\mathit{Nor}(K)\times[-1,1]$ will be important for
describing our limiting processes.
Note, however, that it will eventually be possible to work with the
cylinder $\partial K \times[-1,1]$, which is much easier to visualize.

We also need the so-called local magnification map $\tau_\varepsilon
$; see \cite{K07}. Any point $z \in\mathbb{R}^d \backslash S_K$ can
be written as $z=\Pi(z)+d_s(z)u$, where $d_s(z)$ is the signed (``$+$''
outside) distance between $z$ and $\Pi(z)$ and $u$ an outer normal at
$\Pi(z)$ that satisfies the equality. Now, define
\[
\tau_\varepsilon(z)= \biggl(\Pi(z), u , \frac{d_s(z)}{\varepsilon
} \biggr),\qquad z \in\mathbb{R}^d \backslash S_K , \varepsilon>0.
\]
Observe that $\tau_\varepsilon$ maps $\mathcal{V}_{\varepsilon
}(\partial
K)\backslash S_K$ into $\Sigma$.

We are now prepared to introduce the aforementioned differentiation
of sets. Consider the first support measure $\vartheta_{d-1}$ on
$\mathit{Nor}(K)$; see
\cite{S93}. It attributes measure 0 to the set of all points
$(x,u)$, where, at $x$, there is more than one outer normal $u$. Hence,
we can map it to $\partial K$ in a one-to-one way. On $\partial K$,
this map coincides with Hausdorff measure $\nu$ and if,
for a Borel set $H\subset \mathit{Nor}(K)$, we write
\[
H_0=\{x\in\partial K \dvtx  (x,u) \in H \} ,
\]
\vspace{-4pt}
then
\vspace{-4pt}
%
%
\begin{equation}\label{mes}
\vartheta_{d-1}(H)=\nu(H_0) .
\end{equation}
On $\Sigma$, define the measure $M=\vartheta_{d-1} \times\mu$
($\mu$ being one-dimensional Lebesgue measure). Consider a (Borel)
set-valued function $K(\varepsilon)$, $\varepsilon\in[0,1]$, such
that $K(0)=K$,
with $K$ as before; write
$A(\varepsilon)
=K(\varepsilon)\triangle K$ and assume that $A(\varepsilon)\subset
\mathcal{V}_{\varepsilon
}(\partial K)$. The set-valued function $A(\varepsilon)$, $\varepsilon
\in
[0,1]$, is called \textit{differentiable} at $\partial K$ and
$\varepsilon=0$ if
there exists a Borel set $B \subset\Sigma$ such that $M(\tau
_\varepsilon
A(\varepsilon)\triangle B)\to0$ as $\varepsilon\to0$
(where $\tau_{\varepsilon} A=\{\tau_{\varepsilon}(z)\dvtx  z \in A\}$).
The set $B$ is called the \textit{derivative} of $A(\varepsilon)$ at
$\partial K$.
In this case, we also say that $K(\varepsilon), \varepsilon\in
[0,1]$, is
differentiable with the same derivative and write
\[
\frac{\mathrm{d}}{\mathrm{d}\varepsilon} K(\varepsilon)\bigg|_{\varepsilon=0}=\frac
{\mathrm{d}}{\mathrm{d}\varepsilon} A(\varepsilon)\bigg|_{\varepsilon=0}=B.
\]
Note that $B$ is not unique, but can be changed on a set of $M$-measure 0.

Let $P$ now be as in Section \ref{sec1}. We require that the density $p$ can be
approximated in the neighborhood of $\partial K$ by a function
depending only on $\Pi(z)$ and on whether or not $z\in K$. This latter
possibility is easy to imagine in the change set problems: the limit of
$p(z)$ from inside $K$ can indeed be different from that from outside
if $K$ is the change set. More formally, we require the existence of
two functions, $p_+$ and $p_-$, on $\partial K$ such that, as
$\varepsilon\to0$,
%
%
%
\begin{eqnarray}
\label{p+}\frac{1}{\varepsilon}\int_{\mathcal
{V}_{\varepsilon}
(\partial K) \backslash K}|p(z)-p_+(\Pi(z))|\,\mathrm{d}\mu_d(z)&\to&0,
\\
\label{p-}\frac{1}{\varepsilon}\int_{\mathcal
{V}_{\varepsilon}
(\partial K) \cap K}|p(z)-p_-(\Pi(z))|\,\mathrm{d}\mu_d(z)&\to&0.
\end{eqnarray}
Now, define a measure $M_p$ on $\Sigma$ as follows:
\begin{eqnarray*}
\mathrm{d}M_p(x,u,s)&=&p_+(x)\,\mathrm{d}\vartheta_{d-1}(x,u)\times \mathrm{d}s
\qquad\mbox{for } s>0,
\\
\mathrm{d}M_p(x,u,s)&=&p_-(x)\,\mathrm{d}\vartheta_{d-1}(x,u)\times \mathrm{d}s\qquad
\mbox{for } s\leq0.
\end{eqnarray*}
For convenience, assume that $p_+$ and $p_-$ are bounded (although a
weaker, integrability, condition would suffice). An easy, but
practically interesting, situation occurs when $p_+(x)=c_+$ and
$p_-(x)=c_-$ for all $x \in\partial K$, where $c_+,c_-\geq0$ are two
constants.

The following key result from \cite{K07} shows the
``differentiability of sets in measure'':
if $A(\varepsilon)$ is differentiable at $\partial K$,
then
%
%
\begin{equation} \label{Mp}
\frac{\mathrm{d}}{\mathrm{d}\varepsilon} P(A(\varepsilon
))\bigg|_{\varepsilon
=0}=M_p \biggl(\frac{\mathrm{d}}{\mathrm{d}\varepsilon} A(\varepsilon)\bigg|_{\varepsilon
=0} \biggr).
\end{equation}

\section{Main results}\label{sec4}

Let $\mathcal{A}_{\varepsilon_n}$ be as in Section \ref{sec1} and assume
$M_p(\Sigma
)>0$. Writing $a_n=P(\mathcal{V}_{\varepsilon_n}(\partial K))$, it easily
follows, using (\ref{Mp}), that $a_n/\varepsilon_n\to M_p(\Sigma)$. Hence,
we have, just as for $\varepsilon_n$,
\[
a_n\to0 \quad\mbox{and}\quad na_n\to\infty.
\]
Denote by $\mathcal{B}$ the class of all possible derivatives at
$\varepsilon
=0$ corresponding to $\mathcal{A}_{\varepsilon_n}$, which, by definition,
means that $B \in\mathcal{B}$ if and only if there exists a sequence
of sets $(A_n)^\infty_{n=1}$ with $A_n \in\mathcal{A}_{\varepsilon
_n}$ and
$M(\tau_{\varepsilon_n} A_n\triangle B)\to0$. (Observe that for a thus
converging sequence of Borel subsets of $\Sigma$, the limit set is not
well defined. This limit ``set'' is actually an equivalence class of
sets, defined by the property that for any two sets $B_1, B_2$ in the
class, $M(B_1\triangle B_2)=0$. Out of every such an equivalence class,
we choose one limit (Borel) set $B$, say. Whether or not the conditions
of our results are satisfied will depend on the choices of these $B$'s.
In applications/examples, we should choose natural or appropriate $B$'s
to make the theorems work.)

Consider the local empirical process (\ref{lep}) from Section \ref{sec1}. To
establish its limit in distribution, we need the following steps. Write
$\tau_{\varepsilon}^{-1} C=\{z \in\mathcal{V}_{\varepsilon
}(\partial K)\dvtx  \tau
_{\varepsilon}(z) \in C\}$ for a Borel set $C\subset\Sigma$.
First, using the local magnification map, induce the point process
$\Phi_n$ and the distribution $Q_n$ on $\Sigma\dvtx \Phi_n(C)=\Psi
_n(\tau_{\varepsilon_n}^{-1} C)$ and $Q_n(C)=P_{\varepsilon_n}(\tau
_{\varepsilon_n}^{-1}
C)$. Thus, for any Borel set $C \subset\Sigma$, we can define
%
%
\begin{eqnarray} \label{AB}
v_n(C):\!&=&\frac{1}{\sqrt{na_n}}[\Phi_n(C)-na_nQ_n(C)]\nonumber\\[-8pt]\\[-8pt]
&=&\frac{1}{\sqrt{na_n}}[\Psi_n(\tau_{\varepsilon_n}^{-1}
C)-na_nP_{\varepsilon
_n}(\tau_{\varepsilon_n}^{-1} C)]=z_n(\tau_{\varepsilon_n}^{-1}
C).\nonumber
\end{eqnarray}
Hence, we can define the local empirical processes on two classes of
sets: on $\mathcal{B}_n:=\{\tau_{\varepsilon_n}A\dvtx  A\in\mathcal
{A}_{\varepsilon_n}
\}$, which changes with $n$, and on $\mathcal{B}$, the class of its
limits, or derivative sets, which is fixed. We denote these processes by
\[
v_{n,\mathcal{B}_{n}}:=\{v_n(B)\dvtx  B \in\mathcal{B}_{n}\}
\quad\mbox{and}\quad
v_{n,\mathcal{B}}:=\{v_n(B)\dvtx  B \in\mathcal{B}\} .
\]
(Local empirical processes of the latter type -- i.e. for a fixed
$\mathcal{B}$ -- have been studied in, e.g., \cite{EM97,E97}.
Here, however, our main object is $v_{n,\mathcal
{B}_{n}}$; $v_{n,\mathcal{B}}$ is an auxiliary process, a bridge
between $v_{n,\mathcal{B}_{n}}$ and its limiting process.)
Second, we show that the distribution $Q_n$, which ``governs'' these
processes, converges to the distribution $Q(C)=M_p(C)/M_p(\Sigma)$ and
the processes $v_{n,\mathcal{B}_{n}}$ can be approximated by the
processes $v_{n,\mathcal{B}}$. Next, we verify that $v_{n,\mathcal
{B}}$ converges in distribution to a set-parametric Brownian motion
$W_{\mathcal B}$ and, finally, we note that one can switch from
$W_{\mathcal B}$, given on subsets of $\Sigma$, to its isometric
image, given on the ``easier'' cylinder $\partial K\times[-1,1]$.

Below, we write $C_- = \{(x,u,s) \in C\dvtx  (x,u) \in \mathit{Nor}(K), s \leq0\}$
and $C_+=C\backslash C_-$. Note that in the case $p_+(x)=c_+$ and
$p_-(x)=c_-$ for all $x \in\partial K$, we have
\[
Q(C)=\frac{c_+M(C_+)+c_-M(C_-)}{(c_++c_-)\nu(\partial K)}.
\]
When, for example, $0=c_-<c_+$, we obtain $Q(C)=Q(C_+)=M(C_+)/\nu
(\partial K)$.

For Borel sets $C,C' \subset\Sigma$, define $d(C,C')=(Q(C\triangle
C'))^{1/2}$. Throughout, we will assume that $(\mathcal{B},d)$ is
totally bounded and that
%
%
\begin{equation}\label{1}
\sup_{A \in\mathcal{A}_{\varepsilon_n}} \inf_{B\in\mathcal{B}}
d(\tau
_{\varepsilon_n} A,B)\to0.
\end{equation}
In particular, every sequence $(A_n)^\infty_{n=1}$ with $A_n \in
\mathcal{A}_{\varepsilon_n}$ has a subsequence $(A_{n_k})^\infty
_{k=1}$ such
that for some $B \in\mathcal{B}$, $d(\tau_{\varepsilon_{n_k}}
A_{n_k},B)\to
0$. Assumption (\ref{1}) can be written as
\[
\sup_{B_n \in\mathcal{B}_n} \inf_{B\in\mathcal{B}} d(B_n,B)\to0.
\]
From the definition of $\mathcal{B}$ and the assumption that
$(\mathcal{B},d)$ is totally bounded, it follows that
\[
\sup_{B\in\mathcal{B}} \inf_{B_n \in\mathcal{B}_n} d(B_n,B)\to0.
\]
Thus, the Hausdorff distance between the classes $\mathcal{B}_n$ and
$\mathcal{B}$ tends to 0:
%
%
\begin{equation}\label{gn}\gamma_n:=\max\Bigl(\sup_{B_n \in
\mathcal{B}_n} \inf_{B\in\mathcal{B}} d(B_n,B), \sup_{B\in
\mathcal{B}} \inf_{B_n \in\mathcal{B}_n} d(B_n,B) \Bigr)\to0.
\end{equation}

Recall that it is the aim of this paper to present a central limit
theorem for
$z_{n,\mathcal{A}_{\varepsilon_n}}$, or, equivalently, $v_{n,\mathcal
{B}_{n}}$. By ``central limit theorem for $z_{n,\mathcal
{A}_{\varepsilon
_n}}$'' we mean:
\[
\mbox{(a)}\qquad \sup_{B_n \in\mathcal{B}_n,
B\in\mathcal{B}; d(B_n,B)\leq\gamma_n}|
v_n(B_n)-v_n(B)|\stackrel{P}\to0;
\]
and
\[
\mbox{(b)}\qquad v_{n,\mathcal{B}} \stackrel{d}\to
W_{\mathcal{B}}:=\{W(B), B\in\mathcal{B}\}.
\]
Here, $W_{\mathcal{B}}$ is set-parametric Brownian motion: a bounded,
uniformly $d$-continuous Gaussian process with mean 0 and covariance
structure $\mathbb{E}W(B)W(B')=Q(B\cap B')$. We view $v_n$ and $W$
as processes taking values in $\ell^\infty(\mathcal{B})$ endowed
with the uniform distance and understand weak convergence in the
sense of van der Vaart and Wellner \cite{VW96}.
(We assume, for convenience, that our classes of sets are such that the
various ``suprema'' are measurable, i.e., that they are random variables.)
The following fact is very useful for proving this central limit theorem.
\begin{lemma}\label{lemma1}
From (\ref{p+}) and (\ref{p-}), it follows that $Q_n$ converges to
$Q$ in total variation:
\[
{\sup}|Q_n(C)-Q(C)|\to0,
\]
with the $\sup$ taken over all Borel sets $C\subset\Sigma$.
\end{lemma}

Define $ d_n(A, A'):=({P_{\varepsilon_n}(A \triangle A')})^{1/2}=({P(A
\triangle A')/a_n})^{1/2}$; observe that $d_n(A, A') = (Q_n (\tau
_{\varepsilon_n}A \triangle\tau_{\varepsilon_n}A' ))^{1/2}$.
Assume, for any $\delta>0$, that there exists a finite collection of
pairs (brackets) $[\underline{A}(\delta), \overline{A}(\delta)] $
of Borel sets in $\mathcal{V}_{\varepsilon_n}(\partial K)$ with
$d_n(\underline{A}(\delta),\overline{A}(\delta))\leq\delta$, such that
any set $A \in\mathcal{A}_{\varepsilon_n}$ can be placed in a
bracket from
this collection: $\underline{A}(\delta)\subset A \subset\overline
{A}(\delta) $. Consider such a class of brackets with minimal
cardinality; denote this cardinality (the bracketing number) by $N_{[\cdott],n}(\delta)$ and let $\mathcal{N}_{[\cdott],n}(\delta)$ be the set of
$\underline{A}(\delta)$'s in this class. We assume the same for $\tau^{-1}_{\varepsilon_n} \mathcal{B}:=\{\tau^{-1}_{\varepsilon_n} B\dvtx  B
\in\mathcal{B}\}$ and use the notation $\tilde N_{[\cdott],n}(\delta)$ and
$\tilde{\mathcal{N}}_{[\cdott],n}(\delta)$. We will require
%
%
\begin{eqnarray}
\label{N}\lim_{\delta\downarrow0}\limsup_{n\to
\infty}
\int_0^\delta\sqrt{\log N_{[\cdott],n}(x)} \,\mathrm{d}x&=&0,\\
\label{n}\lim
_{\delta\downarrow0}\limsup_{n\to\infty}
\int_0^\delta\sqrt{\log\tilde N_{[\cdott],n}(x)} \,\mathrm{d}x&=&0.
\end{eqnarray}
\begin{theorem} \label{theo1}
Under the aforementioned assumptions, in particular, the growth
conditions on $\varepsilon_n$, the approximation of $p$ by $p_+$ or
$p_-$ in
(\ref{p+}) and (\ref{p-}), the relation between $\mathcal
{A}_{\varepsilon
_n}$ and $\mathcal{B}$ specified in (\ref{1}) and the entropy
conditions (\ref{N}) and (\ref{n}), the central limit theorem for
$z_{n,\mathcal{A}_{\varepsilon_n}}$ holds, that is, statements \textup{(a)}
and \textup{(b)}
hold true.
\end{theorem}

We also present a version of Theorem \ref{theo1} without assuming bracketing
conditions. To be more precise, we will assume that our classes of sets
near $\partial K$ are Vapnik--Chervonenkis (VC) classes (see, e.g., \cite{VW96},
Section 2.6, for definition and properties).
\begin{theorem} \label{theo2}
Let $\mathcal{A}_{\varepsilon_n}$ be a VC class with index $t_n\leq
t$ for
some $t \in\mathbb{N}$; also, assume that $\mathcal{B}$ is a VC
class. If we assume that $\varepsilon_n \to0$, $n\varepsilon_n\to
\infty$ and (\ref
{p+}), (\ref{p-}) and (\ref{1}), then the central limit theorem for
$z_{n,\mathcal{A}_{\varepsilon_n}}$ holds, that is, statements
\textup{(a)} and \textup{(b)} hold true.
\end{theorem}
\begin{remark}\label{Remark1}
Consider the canonical example of Section \ref{sec1}
and let $\mathcal{K}$ be a VC class. Then $\mathcal{A}$ is also a VC
class, with index $t$, say. Since $\mathcal{A}_{\varepsilon_n}\subset
\mathcal{A}$, the index $t_n$ of $\mathcal{A}_{\varepsilon_n}$ indeed
satisfies $t_n \leq t$.
\end{remark}
\begin{remark}\label{Remark2}
Similar to the discussions in \cite{K07,KW08},
we note that Theorems \ref{theo1} and \ref{theo2}, as
well as the whole construction, can
be carried over to the case where $K$ is a finite union of convex bodies
and, even more easily, to the case where $K$ is closed and bounded and
has a
boundary of positive reach (intuitively, $K$ has a ``smooth'' boundary).
Indeed, the key objects, such as the local magnification map $\tau
_\varepsilon$
(uniquely defined almost everywhere on $\mathbb{R}^d$), the local
Steiner formula, the
notion of derivative sets and Lemma \ref{lemma1}, are all valid for such a $K$.
Moreover, the
existence of the local Steiner formula for a very general $K$ has been
demonstrated in \cite{HLW04}. This offers perspectives for
considering such a general $K$ in the statements of our results.
\end{remark}

The limiting process $W_{\mathcal{B}}$ is defined on subsets of the
cylinder $\Sigma=\mathit{Nor}(K)\times[-1,1]$. This cylinder is not easy to
visualize. However,
since the support measure $\vartheta_{d-1}$ depends on $H$ only
through $H_0$ (cf. (\ref{mes})), we have a similar result for the
measure $Q$. That is, if we write, for a
Borel set $C\subset\Sigma$,
\[
C_0=\{(x,s)\in\partial K \times[-1,1]\dvtx  (x,u,s) \in C \}
\]
and we use the same letter $Q$ for the measure
\[
\mathrm{d}Q(x,s)=\frac{p_{\pm}(x)\,\mathrm{d}\nu(x)\times \mathrm{d}s}{\int_{\partial K}
(p_{+}(x)+p_-(x))\,\mathrm{d}\nu(x)}\qquad \mbox{for } s\gtrless0 ,
\]
which lives on $\partial K\times[-1,1]$, then
\[
Q(C_0)=Q(C).
\]
Therefore, if convenient, we will replace $\Sigma$
by $\Gamma=\partial K \times[-1,1]$ and replace $W_{\mathcal{B}}$
with the
process $W_{\mathcal{B}_0}$ defined on $\mathcal{B}_0=\{B_0\dvtx B \in
\mathcal{B}\}$, a
class of subsets of $\Gamma$. However, we could not do this with
$v_{n,\mathcal{B}_n}$.

Weak convergence in function spaces is important because of its
statistical application, the continuous mapping theorem, which states
that continuous functionals, or statistics,
of the random elements involved inherit the weak convergence. We now
formulate a continuous mapping theorem in our unusual setting, where
the domain of the functions depends on $n$. Let
$\ell^\infty(\mathcal{B}_n)$ and $\ell^\infty(\mathcal{B})$ be the
spaces of bounded functions on $\mathcal{B}_n$ and $\mathcal{B}$,
respectively; let $x_n\in\ell^\infty(\mathcal{B}_n)$, $x \in
\ell^\infty(\mathcal{B})$ and assume that $x$ is $d$-continuous. Also,
assume the functionals $\varphi_n\dvtx  \ell^\infty(\mathcal{B}_n) \to
\mathbb{R}$ and $\varphi\dvtx \ell^\infty(\mathcal{B})\to\mathbb{R}$ are
such that (with $\gamma_n$ as in~(\ref{gn}))
%
%
\begin{equation}\label{xnx}\sup_{B_n \in\mathcal{B}_n, B\in
\mathcal{B}; d(B_n,B)\leq\gamma_n}|x_n(B_n)-x(B)|
\to0
\end{equation}
implies
\[
\varphi_n(x_n)\to\varphi(x).
\]
We then have
%
%
\begin{equation}\label{cmnn}
\varphi_n(v_{n,\mathcal{B}_n}) \stackrel{d}\to
\varphi(W_{\mathcal{B}})
.
\end{equation}
As an example, we see that
\[
\sup_{B_n \in\mathcal{B}_n}|v_n(B_n)|\stackrel{d}\to\sup_{B \in
\mathcal{B}}|W(B)|.
\]

For the proof of (\ref{cmnn}), we only mention that a Skorokhod almost
sure representation theorem yields the existence of $\tilde v_{n,
\mathcal{B}}\stackrel{d}{=}v_{n, \mathcal{B}}$ and $\tilde
W_\mathcal{B}\stackrel{d}{=}W_\mathcal{B}$ such that
\[
\sup_{B\in\mathcal{B}}|\tilde v_n(B)-\tilde W(B)|\to0 \qquad\mbox{a.s.}
\]
If we extend $\tilde v_{n, \mathcal{B}}$ to $\mathcal{B}_n$, we
obtain, from (a),
\begin{eqnarray*}
&&\sup_{B_n \in\mathcal{B}_n, B\in
\mathcal{B}; d(B_n,B)\leq\gamma_n}|\tilde v_n(B_n)-\tilde W(B)| \\
&&\quad\leq \sup_{B_n \in\mathcal{B}_n, B\in\mathcal{B}; d(B_n,B)\leq
\gamma_n}|\tilde v_n(B_n)-\tilde v_n(B)|\\
&&\qquad{}+\sup_{B \in\mathcal{B}}|\tilde v_n(B)-\tilde W(B)|\stackrel{P}\to0.
\end{eqnarray*}
Now, compare this with (\ref{xnx}).
The rest of the proof is elementary.
\begin{Example}\label{Example1}
Let $K=\{(x,y)\in\mathbb{R}^2\dvtx
x^2+y^2\leq1\}$ be the unit disc, so $\partial K=\{(x,y)\in
\mathbb{R}^2\dvtx  x^2+y^2=1\}$ is the unit circle. We have
$S_K=\{(0,0)\}$ and $r(x)=1$ for all $x \in
\partial K$. Also, $\mathcal{V}_{\varepsilon}(\partial K)=\{(x,y)\in
\mathbb{R}^2\dvtx  (1-\varepsilon)^2\leq x^2+y^2\leq(1+\varepsilon)^2\}$.

\textup{(a)} Let $\mathcal{E}$ be the VC class of all closed ellipses (with
interior) in $\mathbb{R}^2$. This $\mathcal{E}$ is an example of the
general $\mathcal{K}$ in the canonical example in Section \ref{sec1}. Thus,
$\mathcal{A}=\{E \triangle K\dvtx  E \in\mathcal{E}\}$ and
$\mathcal{A}_{\varepsilon_n}=\{A\in\mathcal{A}\dvtx  A \subset
\mathcal{V}_{\varepsilon_n}(\partial K)\}$. By Remark \ref{Remark2},
$\mathcal{A}_{\varepsilon_n}$ is a VC class with uniformly bounded index.

We parametrize $\partial K$ with the angle $\theta\in[0, 2\uppi)$
and re-express the cylinder $\Gamma=\partial K\times[-1,1]$ as $[0,
2\uppi)\times[-1,1]$. Consider the functions $f_{\alpha, a,
b,c,d}\dvtx [0, 2\uppi)\to[-1,1]$, defined by
\[
f_{\alpha, a, b,c,d}(\theta)=f(\theta)=a +b\sin^2(\theta-\alpha
)+c\sin(\theta-\alpha)+d\cos(\theta-\alpha),
\]
with $\alpha\in[0,\uppi/2)$ and $a,b,c,d\in\mathbb{R}$ such that
$\sup_{0\leq\theta< 2\uppi} |f_{0, a, b,c,d}(\theta)|\leq1$. Denote
the class of all such functions by $\mathcal{F}_\mathcal{E}$. A
tedious calculation shows that
\[
\mathcal{B}_0=\bigl\{\{(\theta,y)\in[0,2\uppi)\times[-1,1]\dvtx  0< y\leq
f(\theta) \mbox{ or } f(\theta)< y\leq0 \}\dvtx f\in\mathcal
{F}_\mathcal{E} \bigr\}.
\]
Since $\mathcal{B}_0$ is a limit class, it can be shown, directly
using the definition of a VC class, that $\mathcal{B}_0$ is also a VC
class. For $B \in\mathcal{B}_0$, note that for every $\theta\in
[0,2\uppi)$, the intersection of $B$ with $\{(\theta,y)\dvtx y \in[-1,1]\}$
is convex (an interval). Part (b) shows that this need not be the case
in general.

\textup{(b)} Consider, for the same $K$, the very simple class
\[
\mathcal{A}_{\varepsilon_n}=\bigl\{\{ z\in\mathbb{R}^2\dvtx  \Vert z
-\partial K
\Vert/\varepsilon_n\in[a,b]\cup[c,d]\} \dvtx  -1\leq a \leq b \leq c
\leq d
\leq1\bigr\}.
\]
Now,
\[
\mathcal{B}_0=\bigl\{\{(\theta,y)\in[0,2\uppi)\times[-1,1]\dvtx  y \in
[a,b]\cup[c,d] \} \dvtx  -1\leq a \leq b \leq c \leq d \leq1\bigr\}.
\]
Here, $\mathcal{B}_n=\mathcal{B}$.
\end{Example}
\begin{Example}\label{Example2}
Let $K=\{(x,y)\in\mathbb{R}^2\dvtx  0\leq
x,y\leq1\}$ be the unit square with boundary $\partial K$. We
obtain $S_K=\{(x,x)\dvtx 0<x<1\}\cup\{(x,1-x)\dvtx 0<x<1\}$ and for, for example,
$\{(x,0)\dvtx 0\leq x\leq1\}\subset\partial K$, we see that
$r((x,0))=\min(x,1-x)$. It is notationally somewhat cumbersome to
describe $\mathcal{V}_{\varepsilon}(\partial K)$ explicitly, but it is
trivial to see that it is the difference of a set which is a ``square
with circular corners'' and a smaller square.

\textup{(a)} Let $\mathcal{Q}$ be the VC class of all closed quadrangles in
$\mathbb{R}^2$. Set $\mathcal{A}=\{Q \triangle K\dvtx  Q \in
\mathcal{Q}\}$ and $\mathcal{A}_{\varepsilon_n}=\{A\in\mathcal{A}\dvtx  A
\subset\mathcal{V}_{\varepsilon_n}(\partial K)\}$. Again by Remark \ref{Remark2},
$\mathcal{A}_{\varepsilon_n}$ is a VC class with uniformly bounded index.
The present example is somewhat similar to Example \ref{Example1}, but there is a
substantial difference since a square is less smooth than a disc.

We parametrize $\partial K$ with $\theta\in[0, 4)$, the
counterclockwise ``distance'' from the origin, and re-express the
cylinder $\Gamma$ as $[0, 4)\times[-1,1]$. Consider the functions
$f_{a,b}\dvtx [0, 4)\to[-1,1]$, with $a=(a_0,a_1,a_2,a_3)$ and
$b=(b_0,b_1,b_2,b_3)$, defined by
\[
f_{a, b}(\theta)=f(\theta)=a_m(\theta-m)+b_m \qquad \mbox{for }
m \leq\theta< m+1, m=0,1,2,3,\vadjust{\goodbreak}
\]
with $a, b$ such that $a_m\in[-2,2]$ and $\sup_{0\leq\theta< 4}
|f_{a,b}(\theta)|\leq1$. Denote the class of all such functions
by $\mathcal{F}_\mathcal{Q}$. Note that $f \in
\mathcal{F}_\mathcal{Q}$ is typically discontinuous, in contrast to
an $f \in\mathcal{F}_\mathcal{E}$ of Example~\ref{Example1}. It can be shown
that
\[
\mathcal{B}_0=\bigl\{\{(\theta,y)\in[0,4)\times[-1,1]\dvtx  0< y\leq f(\theta
) \mbox{ or }
f(\theta)< y\leq0 \}\dvtx f\in\mathcal{F}_\mathcal{Q} \bigr\}.
\]
It readily follows that $\mathcal{B}_0$ is a VC class.

\textup{(b)} Consider (for the same $K$) a larger class than $\mathcal{Q}$,
namely $\mathcal{C}$, the class of all convex bodies in
$\mathbb{R}^2$. For convenience, let $P$ be the uniform distribution
on $[-1,2]^2$. The class $\mathcal{C}$ is again an example of the
general $\mathcal{K}$ in the canonical example in Section \ref{sec1}, but it is
not a VC class. We have
$\mathcal{A}=\{C \triangle K\dvtx  C \in\mathcal{C}\}$ and
$\mathcal{A}_{\varepsilon_n}=\{A\in\mathcal{A}\dvtx  A \subset
\mathcal{V}_{\varepsilon_n}(\partial K)\}$.

Consider the functions $f \dvtx [0, 4)\to[-1,1]$ defined by
\[
f(\theta)=f_m(\theta-m) \qquad \mbox{for } m \leq\theta< m+1,
m=0,1,2,3,
\]
with $f_m\dvtx  [0,1) \to[-1,1]$ a concave function. Denote the class
of all such functions by $\mathcal{F}_\mathcal{C}$. It can be
shown that
\[
\mathcal{B}_0=\bigl\{\{(\theta,y)\in[0,4)\times[-1,1]\dvtx  0< y\leq f(\theta
) \mbox{ or }
f(\theta)< y\leq0 \}\dvtx f\in\mathcal{F}_\mathcal{C} \bigr\}.
\]
The conditions of Theorem \ref{theo1} are satisfied. In particular, using \cite{VW96},
Corollary 2.7.9, it can be deduced that (\ref{N}) and (\ref{n}) hold true.
\end{Example}

\section{Proofs}\label{sec5}

\begin{pf*}{Proof of Lemma \ref{lemma1}}
Based on the local Steiner formula, in the
proof of Theorem 2 of \cite{KW08},
it is shown that the measure
$P(\tau^{-1}_{\varepsilon_n} \cdot)/\varepsilon_n$
converges in total variation to the measure $M_p$. This implies that
$P(\mathcal{V}_{\varepsilon_n} (\partial K))/\varepsilon_n \to
M_p(\Sigma)$
and hence that
$Q_n=P(\tau^{-1}_{\varepsilon_n} \cdot))/ P(\mathcal{V}_
{\varepsilon_n}(\partial K))$
converges in total variation to $Q=M_p/M_p(\Sigma)$.
\end{pf*}
\begin{pf*}{Proof of Theorem \ref{theo1}} First, we prove statement (a):
\[
\sup_{B_n \in\mathcal{B}_n, B\in\mathcal{B}; d(B_n,B)\leq\gamma
_n}|v_n(B_n)-v_n(B)|\stackrel{P}\to0.
\]
From relation (\ref{AB}), Lemma \ref{lemma1}
and the Markov inequality, it follows that it is sufficient to show that
%
%
\begin{equation}\label{E}\lim_{\delta\downarrow0}\limsup_{n\to
\infty}\mathbb{E}\mathop{\sup_{A \in\mathcal{A}_{\varepsilon_n},
\tilde A \in\tau^{-1}_{\varepsilon_n} \mathcal{B}}}_{ d_n(A, \tilde A)<
\delta}|z_n(A)-z_n(\tilde A)|= 0.
\end{equation}

We use \cite{Vaart98}, Lemma 19.34, page 286, for the proof of
(\ref{E}); in that lemma, we choose the indexing functions to be
$1_A-1_{\tilde A}$.
We then obtain, taking the $\delta$ there to be equal to $\delta\sqrt
{a_n}$, that for some constant $c$,
\begin{eqnarray*}
&&\mathbb{E}\mathop{\sup_{A \in\mathcal
{A}_{\varepsilon
_n}, \tilde A \in\tau^{-1}_{\varepsilon_n} \mathcal{B}}}_{ d_n(A,
\tilde A)<
\delta}|z_n(A)-z_n(\tilde A)| \\
&&\quad\leq c \Biggl(
\frac{1}{\sqrt{a_n}} \int_0^{\delta\sqrt{a_n}} \sqrt{\log
\biggl(N_{[\cdott],n} \biggl(\frac{\varepsilon}{\sqrt{a_n}} \biggr)
\tilde N_{[\cdott],n} \biggl(\frac{\varepsilon}{\sqrt{a_n}}
\biggr) \biggr)}\,
\mathrm{d}\varepsilon\\
&&\qquad\hspace*{10.6pt}{} + \frac{\sqrt n}{\sqrt
{a_n}}\int
_{\mathcal{V}_{\varepsilon_n}(\partial K)}1_{ \{
\sqrt{\log( N_{[\cdott],n} (\delta)
\tilde N_{[\cdott],n} (\delta) )}>\delta\sqrt{na_n}
\}}\,\mathrm{d}P \Biggr).
\end{eqnarray*}
Using $na_n\to\infty$ and (\ref{N}), (\ref{n}), we see that the
second term on the right is equal to 0 for small $\delta$ and large
$n$. The first term is easily seen to be bounded by
\[
c
\int_0^{\delta} \sqrt{\log N_{[\cdott],n} (x )
} \,\mathrm{d}x +c
\int_0^{\delta} \sqrt{\log\tilde N_{[\cdott],n} (x )
} \,\mathrm{d}x.
\]
Hence, (\ref{E}) follows using (\ref{N}) and (\ref{n}).

For a proof of statement (b), we need weak convergence of the
finite-dimensional distributions and tightness of $v_{n, \mathcal
{B}}$. The weak convergence of the finite-dimensional distributions
follows easily from Lemma \ref{lemma1} and an appropriate version of the
multivariate central limit theorem.

To prove tightness, we use \cite{VW96}, Theorem
2.11.9, a general bracketing central limit theorem. We will choose $d$
for the semimetric $\rho$ on $\mathcal{B}$ which is required in that
theorem. For tightness, three conditions have to be fulfilled. The
first one holds trivially since $\Psi_n$ is a sum of indicators. The
third one follows readily since it is essentially our condition (\ref
{n}). It remains to show the second condition:
\[
s_n:=\mathop{\sup_{B,B'\in\mathcal{B}}}_{d(B,B')<\delta_n} \sum
_{i=1}^n \mathbb{E} \biggl(\frac{1}{\sqrt{na_n}}1_{\tau
^{-1}_{\varepsilon
_n}B}(X_i)-\frac{1}{\sqrt{na_n}}1_{\tau^{-1}_{\varepsilon
_n}B'}(X_i)
\biggr)^2\to0\qquad \mbox{for every } \delta_n \downarrow0.
\]
However,
\vspace{-5pt}
\begin{eqnarray*}
s_n&=&\frac{1}{na_n} \sup_{d(B,B')<\delta
_n} \sum_{i=1}^n \mathbb{E} 1_{\tau^{-1}_{\varepsilon_n}B \triangle
\tau
^{-1}_{\varepsilon_n}B'} (X_i)
\\
&=& \frac{1}{a_n} \sup_{d(B,B')<\delta_n} P (\tau
^{-1}_{\varepsilon_n}(B \triangle B') )=
\sup_{Q(B \triangle B')<\delta_n^2} Q_n (B \triangle B' ).
\end{eqnarray*}
Lemma \ref{lemma1} now immediately yields $s_n\to0$.
\end{pf*}
\begin{pf*}{Proof of Theorem \ref{theo2}} Again, we first prove statement
(a) and note that it suffices to show,
for any $\eta>0$, that for $\delta>0$ small enough and then for large $n$,
%
%
\begin{equation}\label{tight2} \mathbb{P}\Bigl(\mathop{\sup_{A \in
\mathcal{A}_{\varepsilon_n}, \tilde A \in\tau^{-1}_{\varepsilon_n}
\mathcal{B}}}_{
d_n(A, \tilde A)\leq\sqrt{\delta}}|z_n(A)-z_n(\tilde A)|>2\eta
\Bigr)\leq2 \eta.
\end{equation}
We have, for $n$ large enough,
%
%
\begin{eqnarray} \label{atoc}
&&\mathbb{P}\Bigl(\mathop{\sup_{A \in\mathcal{A}_{\varepsilon_n},
\tilde A
\in\tau^{-1}_{\varepsilon_n} \mathcal{B}}}_{d_n(A, \tilde A)\leq
\sqrt
{\delta}}|z_n(A)-z_n(\tilde A)|>2\eta\Bigr)\nonumber\\
&&\quad=
\mathbb{P}\Bigl(\mathop{\sup_{A \in\mathcal{A}_{\varepsilon_n},
\tilde A \in
\tau^{-1}_{\varepsilon_n} \mathcal{B}}}_{P(A \triangle\tilde A)\leq
\delta
a_n}|z_n(A)-z_n(\tilde A)|>2\eta\Bigr)\nonumber\\
&&\quad=
\mathbb{P}\Bigl(\mathop{\sup_{A \in\mathcal{A}_{\varepsilon_n},
\tilde A \in
\tau^{-1}_{\varepsilon_n} \mathcal{B}}}_{P(A \triangle\tilde A)\leq
\delta
a_n}|z_n(A\backslash\tilde A)-z_n(\tilde A\backslash A)|>2\eta\Bigr)\nonumber
\\[-8pt]\\[-8pt]
&&\quad\leq
\mathbb{P}\Bigl(\mathop{\sup_{A \in\mathcal{A}_{\varepsilon_n},
\tilde A \in
\tau^{-1}_{\varepsilon_n} \mathcal{B}}}_{P(A \triangle\tilde A)\leq
\delta
a_n}|z_n(A\backslash\tilde A)|>\eta\Bigr)\nonumber
\\
&&\qquad{}+
\mathbb{P}\Bigl(\mathop{\sup_{A \in\mathcal{A}_{\varepsilon_n},
\tilde A \in
\tau^{-1}_{\varepsilon_n} \mathcal{B}}}_{P(A \triangle\tilde A)\leq
\delta
a_n}|z_n(\tilde A\backslash A)|>\eta\Bigr)\nonumber\\
&&\quad
\leq2\mathbb{P}\Bigl(\sup_{C \in\mathcal{C}_{n}, P(C)\leq\delta
a_n}|z_n(C)|>\eta\Bigr),\nonumber
\end{eqnarray}
where $\mathcal{C}_{n}=\{A\backslash\tilde A\dvtx A
\in\mathcal{A}_{\varepsilon_n}, \tilde A \in\tau ^{-1}_{\varepsilon_n}
\mathcal {B}\}\cup \{\tilde A\backslash A\dvtx A
\in\mathcal{A}_{\varepsilon_n}, \tilde A \in\tau ^{-1}_{\varepsilon_n}
\mathcal {B}\}$. It can be shown (see, e.g., \cite{VW96}, page 147),
using $A_1\backslash A_2= A_1 \cap A_2^c$, that $\mathcal{C}_n$ is a VC
class. Also, the index $w_n$ of this VC class is bounded: $\max
_{n\in\mathbb{N}} w_n < \infty$.

We have, writing $N=\Psi_n(\mathcal{V}_{\varepsilon_n}(\partial K))$
and $k=na_n$,
that
\begin{eqnarray*}
&&
\mathbb{P}\Bigl(\sup_{C \in\mathcal{C}_{n}, P(C)\leq\delta
a_n}|z_n(C)|>\eta\Bigr) \\
&&\quad=
\sum_{m=0}^n \mathbb{P}\Bigl(\sup_{C \in\mathcal{C}_{n}, P(C)\leq
\delta a_n}|z_n(C)|>\eta\big| N=m\Bigr) \mathbb{P}(N=m)\\
&&\quad=
\sum_{m=0}^n \mathbb{P} \biggl(\sup_{C \in\mathcal{C}_{n}, P(C)\leq
\delta a_n} \biggl|\frac{1}{\sqrt{k}}[\Psi_n(C)-n P(C)]
\biggr|>\eta\Big| N=m \biggr) \mathbb{P}(N=m)
\\
&&\quad\leq
\sum_{m=\lceil k-C_\eta\sqrt k\rceil}^{m=\lfloor k+C_\eta\sqrt k
\rfloor} \mathbb{P} \biggl(\sup_{C \in\mathcal{C}_{n}, P(C)\leq
\delta a_n} \biggl|\frac{1}{\sqrt{k}}[\Psi_n(C)-n P(C)]
\biggr|>\eta\Big| N=m \biggr) \mathbb{P}(N=m)
\\
&&\qquad{} +\mathbb{P} \bigl(|N-k|\geq C_\eta\sqrt k\bigr),
\end{eqnarray*}
where $C_\eta$ is chosen such that the latter probability concerning
the $\operatorname{binomial}(n, k/n)$ random variable $N$ is bounded by $\eta/2$
for large $n$. Hence, for large $n$,
%
%
\begin{eqnarray} \label{vcbound}
&&
\mathbb{P}\Bigl(\sup_{C \in\mathcal{C}_{n}, P(C)\leq\delta
a_n}|z_n(C)|>\eta\Bigr) \nonumber\\
&&\quad\leq
\sum_{m=\lceil k-C_\eta\sqrt k\rceil}^{m=\lfloor k+C_\eta\sqrt k
\rfloor} \mathbb{P} \Biggl(\sup_{C \in\mathcal{C}_{n},
P_{\varepsilon
_n}(C)\leq\delta} \Biggl|\frac{1}{\sqrt{m}}\Biggl[\sum
_{j=1}^m1_C(Y_j)-mP_{\varepsilon_n}(C)\Biggr] \Biggr|>\frac{\eta}{3}
\Biggr)
\mathbb{P}(N=m)\qquad\\
&&\qquad{}+\sum_{m=\lceil k-C_\eta\sqrt k\rceil}^{m=\lfloor k+C_\eta\sqrt k
\rfloor} \mathbb{P} \biggl(\sup_{C \in\mathcal{C}_{n},
P_{\varepsilon
_n}(C)\leq\delta}\frac{1}{\sqrt{k}}|m-k|P_{\varepsilon_n}(C)>\frac
{\eta
}{2} \biggr) \mathbb{P}(N=m)+\frac{\eta}{2},\nonumber
\end{eqnarray}
where the $Y_j$ are i.i.d. random vectors on $\mathcal
{V}_{\varepsilon
_n}(\partial K)$ distributed according to
$P_{\varepsilon_n}$. Note that in the first probability of the
second sum,
no randomness is involved and that this sum is equal to 0 for $\delta$
small enough. For the first sum, we need a good bound for exceedance
probabilities for the supremum of the empirical process on a VC class.
We will use \cite{A84}, Corollary 2.9. Using $\max_{n \in
\mathbb{N}} w_n <\infty$, this leads to the following upper bound for
the left-hand side of (\ref{vcbound}):
\[
\sum_{m=\lceil k-C_\eta\sqrt k\rceil}^{m=\lfloor k+C_\eta\sqrt k
\rfloor}16\exp\bigl(-\eta^2/(36\delta)\bigr)\mathbb{P}(N=m)+\frac{\eta
}{2}\leq16\exp\bigl(-\eta^2/(36\delta)\bigr)+\frac{\eta}{2}\leq\eta
\]
for small enough $\delta$. So, because of (\ref{atoc}), we have
proven (\ref{tight2}) and hence (a).

For a proof of (b), we only need to show tightness of $v_{n, \mathcal
{B}}$ since the weak convergence of the finite-dimensional
distributions follows as in the proof of Theorem \ref{theo1}.

For proving tightness, we need that, for any $\eta>0$,
\[
\lim_{\delta\downarrow0}\limsup_{n\to\infty}\mathbb{P}\Bigl(\mathop{\sup
_{
B,B' \in\mathcal{B}}}_{ d(B, B')\leq\delta}|v_n(B)-v_n(B')|>\eta\Bigr)= 0
\]
(see, e.g., \cite{VW96}, Theorem 1.5.7). Again,
from (\ref{AB}) and Lemma \ref{lemma1}, it suffices to show that
%
%
\begin{equation}\label{vctight}\lim_{\delta\downarrow0}\limsup
_{n\to\infty}\mathbb{P}\Bigl(\mathop{\sup_{A , A' \in\tau
^{-1}_{\varepsilon
_n} \mathcal{B}}}_{d_n(A, A')\leq2\delta}|z_n(A)-z_n(A')|>\eta\Bigr)= 0.
\end{equation}
The proof of (\ref{vctight}) can be given along the same lines as the
proof of (\ref{tight2}).
\end{pf*}

\section*{Acknowledgements}
We are grateful to two referees for thoughtful comments that led to
improvements of this paper.

\printhistory

\end{document}